\documentclass[11pt]{amsart}
\usepackage{amssymb}
\renewcommand{\a}{\alpha}
\renewcommand{\b}{\beta}
\renewcommand{\d}{\delta}
\newcommand{\g}{\gamma}
\newcommand{\e}{\varepsilon}
\newcommand{\f}{\varepsilon}
\newcommand{\h}{\theta}
\renewcommand{\l}{\lambda}
\newcommand{\s}{\sigma}

\renewcommand{\o}{\omega}
\newcommand{\p}{\phi}
\newcommand{\q}{\psi}

\renewcommand{\O}{\Omega}

\newcommand{\D}{\Delta}

\renewcommand{\P}{\Phi}

\newcommand{\SC}{{\mathcal{C}}}

\newcommand{\SH}{{\mathcal{H}}}

\newcommand{\SO}{{\mathcal{O}}}
\newcommand{\SP}{{\mathcal{P}}}

\newcommand{\PP}{\mathbb{P}}

\newcommand{\C}{\mathbb{C}}

\newcommand{\R}{\mathbb{R}}

\newcommand{\CP}{\mathbb{CP}}

\newcommand{\isom}{\cong}

\newcommand{\codimC}{\text{codim}_{\C}\,}
\newcommand{\surj}{\twoheadrightarrow}

\newcommand{\GL}{\operatorname{GL}}

\newcommand{\sref}{s_{k,x}^{\text{\scriptsize ref}}}
\newcommand{\bd}{\partial}
\newcommand{\bbd}{\bar{\partial}}
\newcommand{\x}{\times}
\newcommand{\ox}{\otimes}

\newtheorem{proposition}{Proposition}[section]
\newtheorem{theorem}[proposition]{Theorem}
\newtheorem{definition}[proposition]{Definition}
\newtheorem{lemma}[proposition]{Lemma}

\newtheorem{remark}[proposition]{Remark}

\title{Codimension one symplectic foliations}

\thanks{Partially supported by The European Contract Human
Potential Programme, Research Training Network HPRN-CT-2000-00101.
Third author has conducted his research with a grant from
Fundaci\'on Pedro Barri\'e de la Maza.}

\date{October, 2001.}

\author{Omegar Calvo}
\address{Cimat: Ap.\ Postal 402 \\ Guanajuato \\
Gto.\ 36000 M\'{e}xico} \email{omegar@fractal.cimat.mx}

\author{Vicente Mu\~noz}
\address{Departamento de Matem\'aticas \\
Universidad Aut\'onoma de Madrid
\\ 28049 Madrid \\ Spain}
\email{vicente.munoz@uam.es}

\author{Francisco Presas}
\address{Mathematics Department. Bldg. 380 \\
Stanford University \\ Stanford \\
CA 90304 \\ Spain}
\email{fpresas@math.stanford.edu}

\begin{document}

\renewcommand{\theenumi}{\roman{enumi}}

\begin{abstract}
 We define the concept of symplectic foliation on a symplectic
 manifold and provide a method of constructing many examples,
 by using asymptotically holomorphic techniques.
\end{abstract}

\maketitle

\section{Introduction} \label{sec:introduction}

During the last three decades there has been an active field of
research related to the study of holomorphic foliations over a
complex manifold \cite{CM82, CL96, GM90, Pa89}. To define a
codimension one holomorphic foliation we need to fix a holomorphic
line bundle $L$ over the manifold $M$. Then we choose a non-zero
holomorphic section $\a \in H^0(T^*M \ox L)$, satisfying the
integrability condition:
 \begin{equation}
 \a \wedge d\a=0. \label{eqn:integrab}
 \end{equation}
There is an equivalence relation given by multiplication of $\a$
by no-where zero holomorphic functions, and a holomorphic
foliation is defined as an equivalence class of such integrable
$1$-forms. In what follows we restrict ourselves to the case where
$M$ is compact, so that the set of foliations is a subset in the
projective space $\PP H^0(T^*M\otimes L)$.

In this work, we aim to generalize this notion to the symplectic
category. We give the following definition

\begin{definition} \label{def:foliation}
 A symplectic foliation $\a$ with normal line bundle $L$
 on a symplectic manifold $(M, \omega)$ is a non-zero element of
 $\SC^{\infty} (T^*_{\C}M \otimes_{\C} L)$ which satisfies the
 integrability condition \eqref{eqn:integrab}. Also we impose that the
 set of singularities, defined as $S_{\a}=\{ x \in M \,\big| \,
 \a(x)=0 \}$, is a finite union of symplectic submanifolds of
 real codimension greater or equal to four and whose intersections are
 transverse and symplectic. Finally, we impose that for any
 $p \in M-S_{\a}$ the subspace $\ker \a(p)\subset T_pM$ is
 symplectic.

 Two symplectic foliations $\a_1$ and $\a_2$ are considered
 equivalent whenever there is an isomorphism $\q:L \to L$ as
 real plane bundles such that $\q^* \a_2= \a_1$.
\end{definition}

To understand $\ker \a(p)$ as a subspace of $T_pM$, we look at the
isomorphism $T^*_{\C}M\ox_{\C} L=T^*M\ox_{\R} L$, where
$T^*_{\C}M$ is the complexified cotangent bundle. Therefore we may
interpret $\a(p):T_pM\to L_p$ as a real linear map and $\ker
\a(p)\subset T_pM$ is a codimension two subspace.

Now if $\a_1$ and $\a_2$ are equivalent then $S_{\a_1}=S_{\a_2}$
and the topological foliations coincide $\ker {\a_1}=\ker {\a_2}$.
Note that the isomorphism $\q:L\to L$ takes values in $\GL(2,\R)$,
so in particular if there is a nowhere zero complex function $f$
such that $\a_1= f\a_2$ then the foliations are equivalent.

The simplest examples of symplectic foliations are given by the
Lefschetz pencils constructed by Donaldson \cite{Do99}. A chart
$\p:U\subset M \to \C^n$ will be called adapted at the point $x\in
U$ if $(\p_*)_x \o = \o_0$, where $\o_0$ is the standard
symplectic form in $\C^n$. A symplectic Lefschetz pencil on a
$2n$-dimensional symplectic manifold $(M,\o)$ consists of a
codimension $4$ symplectic submanifold $N\subset M$ and a map
$f:M-N\to \CP^1$ such that locally around $N$ there are adapted
coordinates $(z_1,\ldots,z_n)$ with values in $\C^n$ where $f$ is
written as $z_2/z_1$. Also $f$ has finitely many isolated critical
points around which there are adapted coordinates where
$f=z_1^2+\cdots +z_n^2$. Finally the fibers of $f$ are symplectic
off their singularities. These belong a special kind of foliations
defined as follows

\begin{definition} \label{def:kupkafoliation}
 A symplectic foliation $\a$ on a $2n$-dimensional symplectic
 manifold $(M,\o)$ is of Kupka type if the singular set $S_\a$ is a
 disjoint union of
   \begin{enumerate}
   \item isolated points where there are adapted charts
   $(z_1,\ldots, z_n)$ such that $\a=z_1dz_1 +\cdots + z_ndz_n$.
   \item codimension $4$ smooth symplectic submanifolds such that
   each point has an adapted chart $(z_1,\ldots, z_n)$ with
   $\a=\eta(z_1,z_2)$ for a $1$-form $\eta$ of $2$
   complex variables with $d\eta(0)\neq 0$ and $\eta^{-1}(0)=\{0\}$.
   \end{enumerate}
\end{definition}

We want to show a general construction of symplectic foliations

\begin{theorem}\label{thm:main}
  Let $(M,\o)$ be a symplectic manifold. Then $M$ admits
  symplectic foliations of Kupka type which are not symplectic
  Lefschetz pencils. Also $M$ admits symplectic foliations not of
  Kupka type.
\end{theorem}

The method of construction is a generalization of the techniques
developed in \cite{MPS00}. The structure of the paper is as
follows. In section \ref{sec:holom_foliations} we give the basic
results of the theory of holomorphic foliations. Section
\ref{sec:holo_theory} reviews the asymptotically holomorphic
theory introduced in \cite{Do96} and used in \cite{Au99,MPS00}.
Next in section \ref{sec:ah-foli} we introduce the notion of
foliation in this category and check that asymptotically
holomorphic foliations with some property of transversality give
symplectic foliations. In the following section we move on to
prove that it is possible to obtain asymptotically holomorphic
foliations by embedding $M$ into the projective space $\CP^d$ and
intersecting the image with a given holomorphic foliation of
$\CP^d$. Finally section \ref{sec:examples} is devoted to give
some examples of foliations constructed with these techniques.

\section{Codimension one holomorphic foliations}
\label{sec:holom_foliations}

In this section we discuss briefly the theory of holomorphic
foliations on a compact connected complex manifold $M$. A
codimension one holomorphic foliation with singularities in $M$ is
an equivalence class of holomorphic $\a\in H^0(M,T^*M\ox L)$,
where $L$ is a holomorphic line bundle and $\a\wedge d\a=0$.

Given a foliation $\a$, we say that $p\in M$ is a regular point if
$\a(p)\neq 0$. Otherwise, we say that $p$ is singular. The set
 $$
 S_\a=\{ p\in M \,|\, \a(p)=0\}
 $$
is the singular set. If this set has components of codimension
$1$, let $D$ be the corresponding divisor. Then there exists a
holomorphic section $f$ of $\SO(D)$ such that $\a/f$ is a
foliation whose singularities are of codimension two or more. So
we can always suppose that $\codimC S_\a\geq 2$.

For a regular point $p\in M$ there exists an open neighborhood
$U\subset M$ of $p$ such that $\a$ may be written as
 $$
  \a= h\, df
 $$
in $U$, where $h$ and $f$ are holomorphic functions in $U$. Such
$f$ is called {\em first integral\/} and $h$ an {\em integrating
factor.\/} The leaves of the foliation in $U$ are the level
surfaces of $f$. Globally, the leaves of the foliation $\a$ are
the leaves of the foliation defined in $M-S_\a$. If $V$ is a
compact hypersurface of $M$ such that $V-V\cap S_\a$ is a leaf, in
general, we have $V\cap S_\a\neq\emptyset$. In any case, by abuse
of language, we will say that $V$ is a compact leaf of the
foliation.

\subsection{Kupka singularities}

In this section, we will consider an important class of
singularities which have stability properties under deformations.

\begin{definition}
The {\em Kupka singular set} of the foliation $\a$ consists of the
points
 $$
 K_\a=\{p\in M \,\big|\, \a(p)=0,\,d\a(p)\neq0\}.
 $$
\end{definition}

For every connected component $K\subset K_\a$, there exists a
holomorphic $1$-form
 $$
 \eta=A(z_1,z_2)\,dz_1+B(z_1,z_2)\,dz_2,
 $$
called the {\em transversal type\/} at $K$, defined on a
neighborhood $V$ of $0\in \C^2$ and vanishing only at $0$, an open
cover $\{U_i\}$ of a neighborhood of $K$ in $M$ and a family of
submersions $\varphi_i:U_i\to \C^2,$ such that
 $$
 \varphi_i^{-1}(0) =K\cap U_i,\quad\text{and}\quad
 \a|_{U_i} =\varphi_i^* \eta.
 $$

A foliation $\alpha$ is of {\em Kupka type\/} if $K_\a$ is compact
and connected.

The main examples of foliations of Kupka type are the following:
Let $L_1$ and $L_2$ be holomorphic line bundles on $M,$ such that
$L_1^{\ox p}=L_2^{\ox q}$, where $p$ and $q$ are relatively prime,
positive integers. Given $f_1$ and $f_2$ holomorphic sections of
the line bundles $L_1$ and $L_2$ respectively, the holomorphic
section
 $$
 \a=pf_1\,df_2-qf_2\,df_1\in H^0(M,T^*M\ox L_1\ox L_2),
 $$
is a foliation. Moreover the leaves of the foliation represented
by $\alpha$, are the fibers of the meromorphic map
$\phi=f_1^p/f_2^q.$ We say that the map $\phi$ is a {\em
meromorphic first integral} of the foliation represented by $\a$.

A {\em branched Lefschetz pencil\/} (a {\em Lefschetz pencil\/} if
$p=q=1$) is a meromorphic map satisfying the following conditions:
\begin{enumerate}
\item The holomorphic line bundles $L_1$ and $L_2$ are positive.
\item The hypersurfaces $\{f_1=0\}$ and $\{f_2=0\}$ are smooth,
and meet transversely along a codimension two submanifold $K$.
\item The subvarieties defined by $\lambda f_1^p - \mu f_2^q = 0$
with $[\lambda:\mu] \in\CP^1,$ are smooth on $M-K,$ except for a
finite set of points, where they have just a non-degenerate
critical point.
\end{enumerate}

These foliations are of Kupka type with $K_\a=\{f_1=f_2=0\}$.

\begin{theorem}[\cite{CL92}]
Let $\a$ be a foliation of Kupka type in $\CP^n$, $n\geq3$. $K_\a$
is a complete intersection if and only if
$\a=pf_1\,df_2-qf_2\,df_1$.
\end{theorem}

For foliations on $\CP^n$, $n\geq 6$, it may be shown that any
foliation of Kupka type is a branched Lefschetz pencil.

For the unbranched case, we have the following construction
involving the fundamental group~\cite{Ca00}. Consider a family
$(E_t,\s_t)$ of projectively flat bundles of rank two with section
such that $(E_0,\s_0)= (L_1\oplus L_1,(f_1,f_2))$ is a Lefschetz
pencil. If $L_1$ is sufficiently ample, we are able to prove that
$H^0(E_t)\neq0$, and then, we consider the foliation
$\s_t^*\SH_t$, where $\SH_t$ denotes the flat structure on the
$\CP^1$-bundle $\PP(E_t)$.

It is an open question whether any foliation of Kupka type with
positive normal bundle and transversal type $z_2\, dz_1- z_1\,
dz_2$ may be described as above.

\subsection{Logarithmic foliations}\label{subsec:log}

A {\em holomorphic integrating factor} of a foliation $\a$ is a
holomorphic section  $\varphi\in H^0(M,L)$ such that the
meromorphic 1-form $\O=\frac{\alpha}{\varphi}$ is closed.

\begin{theorem}
Let $M$ be a projective manifold with $H^1(M;\C)=0$, and let
$\varphi=\varphi_1^{r_1}\cdots\varphi_k^{r_k}$ be an integrating
factor of a foliation $\a$. Then
  $$
  \O =\frac{\alpha}{\varphi}=
  \sum_{i=1}^k\lambda_i\frac{d\varphi_i}{\varphi_i}
  +d\left(\frac{\psi}{\varphi_1^{r_1-1}\cdots\varphi_k^{r_k-1}}\right),
  $$
where $\lambda_i\in\C$ and $\psi$ is a holomorphic section of the
line bundle $\SO\left(\sum_{i=1}^k(r_i-1)\{\varphi_i=0\}\right)$.
\end{theorem}

>From this equation, we have that the hypersurfaces
$D_i=\{\varphi_i=0\}$ are compact leaves of the foliation $\a$.
The residue theorem implies the  relation:
 $$
 \sum_{i=1}^k\lambda_i\cdot\lbrack\{\varphi_i=0\}\rbrack=0\in
 H^2(M;\C).
 $$
The integrating factor is {\em reduced} if $r_i=1$. In this case
 $$
 \alpha=\varphi_1\cdots\varphi_k\left(\sum_{i=1}^k\lambda_i
      \frac{d\varphi_i}{\varphi_i}\right),
 $$
we say that the foliation is {\em logarithmic\/}. The singular set
is the union of $D_i\cap D_j$ for all possible $1\leq i<j\leq k$.
The Kupka set is
 $$
 K_\a=S_\a - \bigcup_{1\leq i<j<t\leq k} \left( D_i\cap D_j\cap D_t\right),
 $$
and it is therefore not compact for $k\geq 3$.

\section{Asymptotically holomorphic theory}
\label{sec:holo_theory}

Let $(M,\o)$ be a symplectic manifold with $[\o]/2\pi \in
H^2(M;\R)$ an integer cohomology class. Such a symplectic manifold
will be called {\em of integer class.\/} Fix an almost complex
structure $J$ compatible with $\o$ and denote $g(u,v)=\o(u,Jv)$
the associated metric. Let $L\to M$ be the hermitian line bundle
with connection whose curvature is $-i\o$. The key for a search of
symplectic objects is to look for objects which are close to be
$J$-holomorphic. The asymptotically holomorphic techniques
introduced by Donaldson \cite{Do96} give a method of construction
of such objects by means of increasing the positivity of the
curvature of the bundles, which is achieved by twisting with
$L^{\ox k}$ for large $k$. Let us introduce the main notations
following \cite{Au99,MPS00}.

\begin{definition}\label{def:ah-sequence}
 A sequence of sections $s_k$ of hermitian bundles $E_k$ with
 connections on $M$ is called asymptotically holomorphic if
 $|\nabla^p s_k|=O(1)$ for $p\geq 0$ and $|\nabla^{p-1} \bbd s_k|
 =O(k^{-1/2})$ for $p\geq 1$. The norms are evaluated with respect
 to the metrics $g_k=kg$.
\end{definition}

\begin{definition}\label{def:trans-sequence}
 A section $s_k$ of the bundle $E_k$ is $\eta$-transverse
 to $0$ if for every $x\in M$ such that
 $|s_k(x)|<\eta$ then $\nabla s_k(x)$ has
 a right inverse $\h_k$ such that $|\h_k|<\eta^{-1}$.
\end{definition}

This means that at a point $x$ close to the zero set of $s_k$ the
differential $\nabla s_k(x):T_xM \to (E_k)_x$ is surjective and
that, in the orthogonal to the kernel, this map multiplies the
length of the vectors at least by $\eta$. This guarantees that
$Z_k=Z(s_k)$ is a submanifold with bounded curvature $R_{Z_k}$ (in
the metric $g_k$) uniformly on $k$, and that $T_xZ_k$ is within
distance $O(k^{-1/2})$ of being a complex subspace of $T_xM$. The
condition of $s_k$ being asymptotically holomorphic implies that
$Z_k$ is symplectic for large $k$.

\begin{definition} \label{def:ah-submfd}
 A sequence of submanifolds $S_k \subset M$ is called
 asymptotically holomorphic if
  $$
  \angle_M(TS_k, JTS_k)=O(k^{-1/2}), ~~ |R_{S_k}|=O(1).
  $$
\end{definition}

The angle $\angle_M$ measures the distance, in the grassmannian,
between two subspaces \cite[definition 3.1]{MPS00}. Thus for $k$
large, any element of a sequence of asymptotically holomorphic
submanifolds is symplectic.

Our objective will be to define and construct asymptotically
holomorphic foliations. This will be done by embedding our
manifold $M$ into a projective space and intersecting it with a
holomorphic foliation in it.

\begin{definition} \label{def:ah-embedding}
 A sequence of embeddings $\p_k:M\to \CP^d$ is $\g$-asymptotically
 holomorphic for some $\g>0$ if it satisfies the following
 conditions:
\begin{enumerate}
  \item $d\p_k:T_xM \to T_{\p_k(x)}X$ has a left
    inverse $\h_k$ of norm less than $\g^{-1}$ at every
    point $x\in M$, i.e.,\ $d\p_k$ multiplies the length of
    vectors at least by $\g$.
  \item $|(\p_k)_*J- J_0|_{g_k}= O(k^{-1/2})$ on the subspace
    $(\p_k)_* T_xM$.
  \item $|\nabla^p \p_k|_{g_k}= O(1)$ and
    $|\nabla^{p-1} \bbd \p_k|_{g_k}= O(k^{-1/2})$,
    for all $p\geq 1$.
\end{enumerate}
\end{definition}

\begin{theorem}[\cite{MPS00}] \label{thm:embeddings}
 Let $(M, \omega)$ be a closed symplectic $2n$-dimensional manifold
 of integer class endowed with a compatible almost complex
 structure $J$ and let $s_k$ be an asymptotically holomorphic
 sequence of sections of the vector bundles $\C^{d+1}\ox L^{\ox k}$
 with $d\geq 2n+1$. Then for any $\a>0$ there exists another
 asymptotically holomorphic sequence $\s_k$ and $\g>0$ such that:
\begin{enumerate}
  \item $|s_k-\s_k|_{C^2,g_k}< \a$.
  \item $\s_k$ is $\g$-projectizable, i.e.,\ $|\s_k|\geq \g$ on all
     of $M$, and for all $k$.
  \item $\PP(\s_k)$ is a $\g$-asymptotically holomorphic sequence of
     embeddings in $\CP^{d}$ for $k$ large enough.
  \item $\p_k^* [\o_{FS}]=k[\o]$, where $\o_{FS}$ is the
     Fubini-Study symplectic form on $\CP^d$.
\end{enumerate}
 Moreover, let us have two asymptotically holomorphic sequences
 $\p_k^0$ and $\p_k^1$ of embeddings in $\CP^{d}$, with
 respect to two compatible almost complex structures. Then for
 $k$ large enough, there exists an isotopy of asymptotically holomorphic
 embeddings $\p_k^t$ connecting $\p_k^0$ and $\p_k^1$.
\end{theorem}

The difference with the result stated in \cite{MPS00} is that here
we use $C^2$-close perturbations, but this makes no real
difference. To be able to intersect the embedded manifold with a
complex submanifold of $\CP^d$ having control of the resulting
submanifold, we need a notion of estimated transversality

\begin{definition} \label{def:openangle}
 Let $N\subset \CP^d$ be a complex smooth submanifold and choose a
 distribution of complex subspaces
 $D_N(y)\subset T_y\CP^d$, in a neighborhood of $N$,
 which extends the tangent distribution to $N$. An embedding $\p_k:M\to
 \CP^d$ is $\s$-transverse to $N$, with $\s>0$
 small enough, if for all $x\in M$,
 $$
  d(\p_k(x),N)< \s \Rightarrow \angle_m((\p_k)_*(T_xM),D_N(\p_k(x)))
  >\s.
 $$
\end{definition}

This angle $\angle_m$ measures the amount of transversality
between two intersecting vector subspaces \cite[definition
3.3]{MPS00}. The condition above assures that the intersection
$\p_k(M)\cap N$ is a sequence of asymptotically holomorphic
submanifolds of $M$ (see \cite[proposition 3.10]{MPS00}).
Moreover, this condition may be achieved.

\begin{theorem}[\cite{MPS00}] \label{good_inter1}
 Let $\phi_k= \PP (s_k)$, where $s_k$ is a $\g$-projectizable
 asymptotically holomorphic sequence of sections of $\C^{d+1} \ox
 L^{\ox k}$, $d\geq 2n+1$, such that $\p_k$ is a
 $\g$-asymptotically holomorphic sequence of embeddings, for some
 $\g>0$. Let $N$ be a complex submanifold in $\CP^{d}$. Then for
 any $\d>0$ there exists an asymptotically holomorphic sequence of
 sections $\s_k$ of $\C^{d+1} \ox  L^{\ox k}$ such that
\begin{enumerate}
  \item $|\s_k-s_k|_{C^2,g_k}<\d$.
  \item $\psi_k=\PP (\s_k)$ is a $\eta$-asymptotically holomorphic embedding
    in $\CP^{d}$ which is $\f$-transverse to $N$, for some $\eta>0$ and
    $\f>0$, for $k$ large enough.
\end{enumerate}
 Moreover, let us have two asymptotically holomorphic sequences
 $\p_k^0$ and $\p_k^1$ of embeddings in $\CP^{d}$, with respect to
 two compatible almost complex structures, which are
 $\f$-transverse to $N$. Then for $k$ large enough, there exists an
 isotopy of asymptotically holomorphic embeddings $\p_k^t$ which
 are $\f'$-transverse to $N$, connecting $\p_k^0$ and $\p_k^1$.
\end{theorem}

\section{Asymptotically holomorphic foliations}
\label{sec:ah-foli}

Consider a foliation $\a$ of $M$ as an application $\a:TM \to L$.
The line bundle $L$ is called the normal bundle of the foliation.
Using the almost complex structure on $TM$, we decompose $\a$ in
complex linear and complex anti-linear parts,
 $$
  \a =\a_{1,0} +\a_{0,1}.
 $$
When $\a_{0,1}(x)=0$ the subspace $\ker \a(x)\subset T_xM$ is
complex. Still when $|\a_{0,1}(x)|<|\a_{1,0}(x)|$ the subspace
$\ker \a(x)$ is symplectic.

\begin{definition}\label{def:ah-fol}
 A sequence of foliations $\a_k$ with hermitian normal bundles
 $E_k$ is asymptotically holomorphic if
 $$
  |\nabla^p \a_k|=O(1), \quad |\nabla^p (\a_k)_{0,1}|=O(k^{-1/2}),
  \qquad  p=0,1,2,
 $$
 $$
  |\nabla^{p-1} \bbd (\a_k)_{1,0}|=O(k^{-1/2}),\qquad p=1,2.
 $$
\end{definition}

Also we need a measure of transversality for foliations. This is
provided by the following definition. As in the holomorphic case,
there is a subset of the singular case, $K_{\a_k} \subset
S_{\a_k}$ which is easily controlled, and which in some cases it
reduces to the Kupka set of $\a_k$.

\begin{definition}\label{def:e-regular}
 Let $\g,\e>0$. A sequence of foliations $\a_k$ with hermitian normal bundles
 $E_k$ is $(\g,\e)$-regular if there is a subset
 $K_{\a_k}$ of $S_{\a_k}$ such that
  \begin{enumerate}
  \item $K_{\a_k}$ is a union of (closed) asymptotically holomorphic
  submanifolds whose intersections are transverse and asymptotically
  holomorphic submanifolds.
  \item Let $B_\g^k$ be the tubular neighborhood of radius $\g$
  of $K_{\a_k}$ in $g_k$-norm. Then $\a_k$ defines a regular
  foliation in $B_{\g}^k-K_{\a_k}$ such that for any point $x\in
  B_{\g}^k-K_{\a_k}$ it is satisfied that
  $$
  \angle_M( \ker \a_k, J\ker \a_k) =O(k^{-1/2}), 
  $$
  i.e., the leaves in $B_{\g}^k-K_{\a_k}$ are asymptotically holomorphic.
  \item $(\a_k)_{1,0}$ is $\e$-transverse to zero as a section
  of $T^{1,0}M^*\ox E_k$ over $M-B_{\g}^k$.
  \item For any point $x\in M-B_{\g}^k$, there is a uniform number
  $r>0$ such that in the ball of $g_k$-radius $r$ centered at $x$,
  the foliation can be written as $\a_k= h_k \, df_k$, for some
  trivialization of $E_k$, where $h_k$ and $f_k$ are
  asymptotically holomorphic and $h_k$ is bounded above and below by
  a uniform constant (independent of $k$ and $x$).
  \end{enumerate}
\end{definition}

\begin{proposition} \label{prop:holo_sim}
 Let $\a_k$ be a sequence of asymptotically holomorphic foliations
 with hermitian normal bundles $E_k$
 which are $(\g,\e)$-regular. Then there exists an
 arbitrarily small $C^1$-perturbation of $\a_k$ which is a
 symplectic foliation for $k$ large enough.
\end{proposition}

\begin{proof}
 For $k$ large enough, $\a_k$ defines a symplectic foliation in
 $B_\g^k$ as a simple corollary of the definition (see \cite[section 3.2]{MPS00}).
 In the complementary $M-B_\g^k$ we study the set of bad
 points
 $$
 B(\a_k)=\{x\in M-B_\g^k \,\,\big| \,\, |(\a_k)_{1,0}(x)|\leq |(\a_k)_{0,1}(x)|\}.
 $$
 Notice that this is the set of points where the distribution is
 singular or is not symplectic. We will say that $x\in M-B_\g^k$ is
 a critical point of $\a_k$ if $(\a_k)_{1,0}(x)=0$. We want to
 modify the foliation so that $B(\a_k)$ only consists of finitely
 many isolated critical points in which the foliation has a
 standard model of the form $z_1dz_1+\ldots +z_ndz_n$.

 Using lemma \ref{lem:shape} below we can take $k$ large enough so that
 the set of bad points is included in a disjoint collection of
 balls of uniform size centered at the critical points.
 Then we perturb the foliation in a small
 neighborhood of the set of critical points,
 by using proposition \ref{prop:new}, to obtain a new
 foliation with $B(\a_k)$ equal to the set of critical points and
 such that in a neighborhood of a critical point it has the form
 $\a_k=\sum z_idz_i$.
\end{proof}

\begin{lemma} \label{lem:shape}
 For $k$ large enough the set $B(\a_k)$ is contained in a finite
 set of disjoint balls $B(x_j, c)$ of uniform $g_k$-radius $c$ around the
 critical points $x_j$, such that the balls $B(x_j,2c)$ are
 disjoint and contained in $M-B_{\g/2}^k$.
\end{lemma}

\begin{proof}
Let us see first that the minimal $g_k$-distance between two
critical points of $\a_k$ in $M-B_\g^k$ is bounded below by a
uniform constant $c_0>0$. Suppose that $x_1$ is a critical point,
i.e.,\ $(\a_k)_{1,0}(x_1)=0$. In a neighborhood of uniform radius,
$|(\a_k)_{1,0}(x)|<\e$. The $\e$-transversality implies that
$\nabla (\a_k)_{1,0}(x):TM\to T^*M\ox E_k$ has an inverse. If
there is another critical point $x_2$ in such a ball, then there
is a point $y$ in the segment joining $x_1$ and $x_2$ such that
$\nabla(\a_k)_{1,0}(y) (x_2-x_1)=0$, which is a contradiction.

Given a point $x\in B(\a_k)$ then the distance of $x$ to the set
of critical points is bounded above, for $k$ large enough, by $c$,
where $c>0$ is an arbitrarily small uniform constant. For $k$
large enough $|(\a_k)_{1,0}(x)| \leq |(\a_k)_{0,1}(x)| \leq
Ck^{-1/2}$. By the $\e$-transversality, $\nabla (\a_k)_{1,0}(x):TM
\to T^*M\ox E_k$ has an inverse of norm bounded by $\e^{-1}$. The
bounds in the second derivatives of $\a_k$ allows us to control
the radius where the inverse function theorem applies (see
\cite[lemma 8]{Do99}). Therefore there must be a zero of
$(\a_k)_{1,0}$ in a neighborhood of some uniform radius $c$.
Moreover this $c$ may be made as small as we please just by
increasing $k$.
\end{proof}

\begin{proposition} \label{prop:new}
 Let $\a_k$ be a sequence of asymptotically holomorphic foliations
 which is $(\g,\e)$-regular. Let $x_j\in M-B_\g^k$ be a
 critical point of $\a_k$ such that there are asymptotically
 holomorphic functions $h_k$ and $f_k$ with $h_k$ bounded below and
 $\a_k=h_k\, df_k$ on a neighborhood of uniform $g_k$-radius
 $B(x_j,2c)$. Then there exists an arbitrarily small
 $C^1$-perturbation of $\a_k$ supported in $B(x_j,2c)$, which is
 a symplectic foliation in the annulus $B(x_j,2c)-B(x_j,c)$ and of the
 form $z_1dz_1+\cdots + z_ndz_n$ in $B(x_j,c)$.
\end{proposition}

\begin{proof}
 Given a critical point $x_j\in M-B_\g^k$, i.e.,
 $(\a_k)_{1,0}(x_j)=0$, write $\a_k=h_k\,df_k$ as in the statement.
 Recall that from the proof of lemma \ref{lem:shape}, the constant
 $c$ can be chosen arbitrarily small, just by increasing the first
 $k$ satisfying the property. This implies that $\nabla
 \alpha_k(p)$ is ``approximately constant'' in that ball (by radial
 parallel transport). We can use an asymptotically holomorphic
 chart $\P_k: B_{\C^n}(0,1) \to B_{g_k}(x_j,2c)$ provided in
 \cite[lemma 3]{Au99} to trivialize the manifold at a neighborhood
 of $x_j$. Consider $f_k$, $h_k$ then as functions of $z_1, \ldots,
 z_n$. We may suppose without loss of generality that $f_k(x_j)=0$
 and that $h_k(x_j)=1$.

 Since $(\a_k)_{1,0}$ is $\e$-transverse at $x_j$ and $\a_k$ is
 asymptotically holomorphic, then for $k$ large enough
 $$
  \bd (\a_k)_{1,0}(0): TM \to T^*M\ox E_k
 $$
 has an inverse of norm bounded by $(\e')^{-1}$, i.e.,\ it multiplies
 the length of vectors by an amount at least $\e'$, for some $\e'$
 slightly smaller than $\e$. Since $\bd
 (\a_k)_{1,0} (0)=\bd\bd f_k(0)$ is the complex Hessian of $f_k$,
 we define
 $$
  H= \frac12 \sum \frac{\bd^2 f_k}{\bd z_i\bd z_j}(0) z_iz_j.
 $$
 Consider the following foliation in $B(x_j,2c)$,
 $$
 \tilde{\a}_k=h_k\,dH.
 $$
 Here $dH$ is a holomorphic foliation with respect to the standard
 complex structure $J_0$ on the ball. Since this asymptotically
 holomorphic chart satisfies that the Nihenjius tensor has norm
 $O(k^{-1/2}|z|)$, we have that
 $(\tilde{\a}_k)_{0,1}=O(k^{-1/2}|z|)$. We have that
 $f_k(0)=H(0)=0$, $\bd f_k(0)=dH(0)=0$ and $\bd\bd
 f_k(0)=\nabla\nabla H(0)$. Since both are asymptotically
 holomorphic we have that $|f_k-H|=O(|z|^3+k^{-1/2}|z|)$.
 Analogously $|\a_k-\tilde{\a}_k|=O(|z|^2+k^{-1/2})$.

 Now let $\b$ be a bump function such that $\b(x)=1$ for $x\in
 B(x_j,c)$, $\b(x)=0$ for $x\not\in B(x_j,3c/2)$ and
 $|\nabla \b|=O(c^{-1})$. Define the foliation in the whole of
 $M$,
 $$
  \hat{\a}_k = h_k\, d\left(\b f_k+ (1-\b) H\right)= \b \a_k+(1-\b)
  \tilde{\a}_k +\nabla \b (f_k-H).
 $$
 We want to prove that outside $x_j$,
\begin{equation}
 |(\hat{\alpha}_k)_{1,0}| > |(\hat{\alpha}_k)_{0,1}|. \label{eqn:key}
\end{equation}
 In $B(x_j,c)$, $\hat{\a}_k=\tilde{\a}_k =h_k\,dH=
 h_k (\sum a_{ij}z_idz_j)$. The $\e$-transversality of
 $(\a_k)_{1,0}$ implies that all the eigenvalues of the symmetric
 matrix $(a_{ij})$ have norm bigger than $\e'$. Therefore
 $|\tilde{\a}_k|\geq \e'|z|/2$ (for $c$ small). Also
 $|(\tilde{\a}_k)_{0,1}|\leq C k^{-1/2}|z|$, for some
 constant $C$, so \eqref{eqn:key}
 holds in $B(x_j,c)$. On $B(x_j,3c/2)-B(x_j,c)$ we have
 $$
 |\hat{\a}_k-\tilde{\a}_k| \leq |\a_k-\tilde{\a}_k| + |\nabla \b|
 |f_k-H|=O(c^2+k^{-1/2}).
 $$
In particular, $|(\hat{\a}_k)_{0,1}|=O(c^2+k^{-1/2})$ and
$|(\hat{\a}_k)_{1,0}|\geq \e' c/2 - O(c^2+k^{-1/2})$. Taking $c$
small (but uniformly on $k$) and then $k$ large enough we have
\eqref{eqn:key}.

Finally take $\tilde{h}_k$ to be equal to $h_k$ in
$B(x_j,2c)-B(x_j,3c/2)$ and equal to $1$ in $B(x_j,c)$. Then
 $$
  \tilde{h}_k\, d\left(\b f_k+ (1-\b) H\right)
 $$
satisfies the required properties and in $B(x_j,c)$ it is of the
form $\sum a_{ij}z_idz_j$. A suitable orthonormal change of
coordinates transforms this into $\sum z_idz_i$.
\end{proof}

\begin{remark}
 Note that the perturbed foliation in the proof above is of the
 form $\a_k=h_k\, df_k$ in $B(x_j,2c)$, so the integral
 submanifolds are the level sets $f_k=\l$. In a small neighborhood
 of the singularity, the leaves of the foliation are of the form
 $\sum z^2_i=\l$.
\end{remark}

\section{Construction of asymptotically holomorphic foliations}
\label{sec:existence}

Once introduced all the asymptotically holomorphic machinery, we
are ready to perform our main construction of asymptotically
holomorphic foliations, generalizing the ideas contained in
\cite{MPS00}. Let $(M,\o,J)$ be a $2n$-dimensional symplectic
manifold of integer class with a fixed compatible almost complex
structure. Let $L\to M$ be the hermitian line bundle with
connection whose curvature is $-i\o$.

Take any asymptotically holomorphic sequence of sections
$s_k=(s_k^0,\ldots, s_k^{d})$ of $\C^{d+1}  \ox L^{\ox k}$ such
that $\p_k=\PP(s_k):M\to \CP^{d}$ is a sequence of asymptotically
holomorphic embeddings with $\p_k^* \SO(1)=L^{\ox k}$, whose
existence is guaranteed by theorem \ref{thm:embeddings}.

Now fix a holomorphic foliation $\a\in H^0(T^*\CP^{d} \ox \SO(N))$
in $\CP^{d}$, such that the singular set $S_\a$ is a union of
smooth complex submanifolds intersecting transversely. There are
many examples of such foliations \cite{Ca94}. We want to study the
restriction $\a_k=\p_k^*\a$ of the foliation $\a$ to the sequences
of embeddings $\p_k$ and to prove that for suitable choice of
embeddings we get asymptotically holomorphic foliations which are
$(\g,\e)$-regular.

\begin{proposition} \label{prop:Kupka-trans}
 Let $\p_k$ be an asymptotically sequence of embeddings of $M$ into
 $\CP^{d}$ and let $\a$ be a foliation as above in the projective
 space. Then there is a $C^2$-close sequence of embeddings $\psi_k$
 which is $\g$-transverse to every submanifold
 of the singular set $S_\a$. Then the induced foliation
 $\a_k=\psi_k^*\a$ in $M$ is an asymptotically
 holomorphic foliation and satisfies conditions (i) and (ii) of
 definition \ref{def:e-regular} with $K_{\a_k}=
 \psi_k^{-1}(\psi_k(M)\cap S_\a)$.
\end{proposition}

\begin{proof}
 Let $\p_k$ be an asymptotically sequence of embeddings of $M$ into
 $\CP^{d}$. Write $S_\a=\cup S_i$, where $S_i\subset \CP^d$ are
 smooth complex submanifolds of $\CP^d$. We may apply theorem
 \ref{good_inter1} to perturb $\p_k$ to a $C^2$-close sequence of
 embeddings $\q_k$ which is $\g$-transverse to every submanifold
 $S_i$, for some uniform $\g>0$. This implies that the
 submanifold $\q_k(M)$ intersects $S_i$ along an asymptotically
 holomorphic submanifold by \cite[proposition 3.10]{MPS00}.

 Let $\a_k=\q_k^*\a$ be the induced foliation in $M$. Then $\a_k$
 is asymptotically holomorphic using the asymptotically holomorphic
 bounds of $\q_k$ and the holomorphicity of $\a$. Note that
 $(\a_k)_{1,0}=\a\circ \bd \q_k$ and $(\a_k)_{0,1}=\a\circ \bbd
 \q_k$.

 Now $K_{\a_k}=\q_k^{-1}(\q_k(M) \bigcap S_{\a})$ is a finite
 union of asymptotically holomorphic submanifolds and it is
 included in $S_{\a_k}$. The $\g$-transversality to
 $S_i$ implies that for the points in a neighborhood of radius $\g$
 of $S_i$, the angle between the tangent space of $\q_k(M)$ and the
 distribution $D_{S_i}$ determined by $S_i$ is bigger than $\g$. Now the
 regular leaves of the foliation $\a$ around $S_i$ contain $S_i$
 in its closure, so that one may assume that $D_{S_i}(x) \subset \ker
 \a(x)$ in a neighborhood of $S_i$. We use a linear
 algebra result \cite[proposition 3.5]{MPS00} that says that for
 $U,V,W$ subspaces of a finite dimensional euclidean vector space
 with $V\subset W$ it is satisfied that $\angle_m(U,V) \leq
 \angle_m(U,W)$. Therefore
 $$
  \angle_m(\ker \a(x), T_x\p_k(M))\geq \angle_m(D_{S_i}(x),
  T_x\psi_k(M))\geq \g,
 $$
 for any $x\in \psi^{-1}(B_\g(S_\a))$. This implies that the leaves are
 asymptotically holomorphic in some $B_{c_o\g}$, for a constant
 $c_o>0$. This gives the sought property (maybe after multiplying
 $\g$ by a suitable uniform constant).
\end{proof}

\begin{theorem} \label{thm:main2}
 Let $\p_k$ be a sequence of asymptotically holomorphic
 embeddings of $M$ into $\CP^{d}$. Fix a holomorphic foliation
 $\a\in H^0(T^*\CP^d\ox \SO(N))$ in $\CP^{d}$ as above.
 Then there exists an arbitrarily
 $C^2$-close sequence of embeddings $\q_k$ such that $\a_k=\q_k^*
 \a$ is an asymptotically holomorphic sequence of foliations of
 $M$ with normal bundle $L^{\ox Nk}$,
 which is $(\g,\e)$-regular for uniform $\g,\e>0$.

 Moreover any two such embeddings $\q_k^i$, $i=0,1$, induce
 isotopic foliations $\a^i_k$, for $k$ large enough.
\end{theorem}

\begin{proof}
Recall that $\p_k=\PP(s_k)$ for a $\g$-asymptotically holomorphic
sequence of sections $s_k$ of $L^{\ox k} \ox \C^{d+1}$ which is
$\g$-projectizable. The property of $\p_k$ being
$\g$-asymptotically holomorphic is open in $C^1$-sense, so any
small perturbation will still be $\g/2$-asymptotically
holomorphic. Using proposition \ref{prop:Kupka-trans} we may
assume that $\p_k$ is already $\g$-transverse to (every
submanifold in) $S_\a$ (reducing $\g$ if necessary). The property
of an asymptotically holomorphic embedding being $\g$-transverse
to $S_\a$ is open in $C^1$-sense \cite[definition 3.11]{MPS00}, so
any small perturbation will still be $\g/2$-transverse. Denote by
$B_{\g}^k$ the tubular neighborhood of radius $\g$ of
$K_{\a_k}=\p_k^{-1}(\p_k(M) \bigcap S_{\a})$ in $M$. We need to
perturb $\p_k$ to a sequence of embeddings such that
$(\a_k)_{1,0}$ is $\e$-transverse to zero in $M-B_{\g/2}^k$.

We define the following property for sequences of sections $s_k$
which are $\g/2$-projectizable and such that $\p_k=\PP(s_k)$ is
$\g/2$-asymptotically holomorphic and $\g/2$-transverse to $S_\a$:
$s_k$ satisfies the property $\SP(\e,x)$ if $(\p^*_k\a)_{1,0}$ is
$\e$-transverse as a section of $T^*M\ox L^{\ox Nk}$ at the point
$x$ or else $x\in B_{\g/2}^k$. This property is local and open in
$C^2$-sense (for $\e$ small). 

We want to use the globalization lemma in \cite[proposition
3]{Au99} which states the following: Let $s_k$ be asymptotically
holomorphic sections of $E_k=L^{\ox k}\ox \C^{d+1}$. If we can
obtain for any point $x\in M$ and any $\d>0$ an asymptotically
holomorphic sequence of sections $\tau_{k,x}$ with Gaussian decay
away from $x$ in $C^r$-norm and $|\tau_{k,x}|_{C^r, g_k}<\d$ such
that $s_k+\tau_{k,x}$ satisfies the property $\SP(\s, y)$ for all
$y$ in a ball of uniform radius $B_{g_k}(x,c)$, with $\s= c'\d
(\log( \d^{-1}))^{-p}$, with $c,c',p$ independent of $k$, then,
given any $\d>0$, there exist, for all large enough $k$,
asymptotically holomorphic sections $\s_k$ of $E_k$ such that
$|s_k-\s_k|_{C^r, g_k}< \d$ and the sections $\s_k$ satisfy
$\SP(\eta, x)$ for all $x\in M$ with $\eta>0$ independent of $k$.

The transversality of $\p_k$ to $S_\a$ implies that
$\p_k(M-B^k_{\g/2}) \subset \CP^d-B_{c_o\g}(S_\a)$, for some
uniform constant $c_o>0$. Now fix a finite covering $U_j$ of
$\CP^d - B_{c_o\g}(S_{\a})$ such that in each of the sets $U_j$
one may write $\a=h_j\,df_j$ where $h_j$ is a (holomorphic)
integrating factor and $f_j$ is a first integral.

Let $x\in M- B^k_{\g/2}$. We may choose $c$ small enough so that
$\p_k(B_{g_k}(x,c)) \subset U_j$ for some $j$ (since $|\nabla
\p_k|_{g_k} \leq C$). Also any small perturbation will still be
inside the same open set. Define $f_k^j=f_j\circ \p_k$ and
$h_k^j=h_j\circ \p_k$. Both are asymptotically holomorphic in the
ball. Moreover $\a_k=h_k^j\,df_k^j$. The functions $h_k^j$ are
bounded above and below by fixed constants. Therefore checking
transversality for $(\a_k)_{1,0}$ is equivalent to checking
transversality for $\bd f_k^j$.

With a transformation of $U(d+1)$ in $\C^{d+1}$ we may suppose
that $s_k(x)=(s_k^0(x),0, \ldots, 0)$. As $s_k$ is
$\g$-projectizable and asymptotically holomorphic, we suppose that
$|s_k^0|\geq \g/2$ on $B_{g_k}(x, c)$ (maybe reducing $c>0$). By
\cite[lemma 2]{Au99} there are asymptotically holomorphic sections
$\sref$ of $L^{\ox k}$ with Gaussian decay away from $x$ and with
$|\sref|\geq c_1$ on $B_{g_k}(x,c)$, for some uniform $c_1>0$.

  We use the standard chart $\Psi_0$ in $\CP^{d}$ around
  $p=\p_k(x)=[1:0\cdots:0]$. With respect to this trivialization the map
  $\p_k$ is given locally as
 \begin{eqnarray*}
  \Psi_0\circ \p_k: B_{g_k}(x,c) & \to & \C^{d} \\
  y & \to & \left( \frac{s_k^1(y)}{s_k^0(y)}, \ldots, \frac{s_k^{d}(y)}{s_k^0(y)} \right).
 \end{eqnarray*}

Now we can suppose that $|\bd f_j(0)|>c_2$, for a universal
constant $c_2$ since we are well away from $S_\a$. Also we may
suppose that $\bd f_j(0)=(0,0,\ldots,0, \frac{\bd f_j}{\bd
w_d}(0))$. Therefore $\frac{\bd f_j}{\bd w_d}$ is big enough in a
small neighborhood.

We trivialize $M$ at a neighborhood of $x$ by using the
asymptotically holomorphic charts $\P_k: B_{\C^n}(0,1) \to
B_{g_k}(x,c)$ provided in \cite[lemma 3]{Au99}. We denote by
$f_k^j$ and $h_k^j$ again the corresponding functions defined in a
ball of $\C^n$, which are asymptotically holomorphic. We define
the ``approximately orthogonal basis'' as in \cite{Au99}
 \begin{equation}
  \mu_k^i= \bd \left(z_i\frac{\sref}{s_k^0}\right). \label{eqn:basis}
 \end{equation}
At $x$ it is an orthogonal basis and all the forms are
asymptotically holomorphic. We can use \eqref{eqn:basis} to
locally trivialize the cotangent bundle. In particular we may
write
 $$
  \frac{\bd f_k^j}{(\bd f_j/\bd w_d) \circ\p_k} =
 t_1 \mu_k^1 + \cdots t_n\mu_k^n.
 $$
It is easy to check that $t:B_{\C^n}(0,1) \to \C^n$ defined by
$t=(t_1, \ldots, t_n)$ is asymptotically holomorphic. This is
because $\{ \mu_k^1, \cdots, \mu_k^n \}$ is close to be an
orthogonal matrix, and so all the eigenvalues are bounded below
and above by positive uniform constants. Moreover the amount of
transversality of $t$ and of $\bd f_k^j$ are related by non-zero
uniform constants. So we need only get transversality for $t$.

The main local result is Donaldson's theorem 12 in \cite{Do99}
stating that there exists $w =(w_1,\ldots,w_n) \in \C^n$ with
$|w|<\d$ such that $t-w$ is $\s$-transverse to $0$ over the ball
$B_{\C^n}(0,\frac{9}{10})$, with $\s=\d(\log (\d^{-1}))^{-p}$, for
a universal $p>0$ and $k$ large enough.

Now define the perturbation
  $$
  \tau_{k,x}=(0,\ldots, 0, -\sum w_iz_i\sref).
  $$
This is asymptotically holomorphic, with Gaussian decay away from
$x$ and norm less than $\d$. The asymptotically holomorphic
sequence $\hat{s}_k= s_k+\tau_{k,x}$ has the corresponding
$\hat{t}=t-w$, so the induced $\hat{\a}_k=\hat{\p}^*_k\a$
satisfies that $(\hat{\a}_k)_{1,0}$ is $C'\s$-transverse over
$B_{g_k}(x,c)$, where $C'$ is again another uniform constant. This
concludes the proof.

For the one-parameter case, let $\q^i_{k}:M\to \CP^d$, $i=0,1$, be
two asymptotically holomorphic sequences of embeddings with
respect to two compatible almost complex structures $J_i$, which
are $(\g,\e)$-regular, for some $\g,\e>0$. Consider a
one-parameter family of compatible almost complex structures
$J_t$, $t\in [0,1]$, interpolating between $J_0$ and $J_1$ and let
$s_k^t$ be $J_t$-asymptotically holomorphic sections of $L^{\ox
k}\ox \C^{d+1}$ such that $\q_k^0=\PP(s_k^0)$ and
$\q_k^1=\PP(s_k^1)$.

We initially perturb $s_k^t$ using theorem \ref{good_inter1} so
that all $\q_k^t=\PP(s_k^t)$ are asymptotically holomorphic
embeddings which are $\g$-transverse to $S_\a$ (reducing $\g>0$ if
necessary). Then the argument above works for one-parameter
families of sections depending on $t\in [0,1]$ since all the
ingredients do (see \cite{Au99,Do99,MPS00}). This means that for a
given $\d>0$, there exists, for large enough $k$,
$J_t$-asymptotically holomorphic sections $\s_k^t$ of $L^{\ox
k}\ox \C^{d+1}$ such that $|s_k^t -\s_k^t| <\d$ and $\PP(\s_k^t)$
are $(\g/2,\eta)$-regular asymptotically holomorphic embeddings in
$\CP^d$ for some uniform $\eta>0$. Taking $\d>0$ very small, the
linear segment $u s_k^0 + (1-u) \s_k^0$, $u\in [0,1]$, consists of
sections inducing $(\g/2,\e/2)$-regular maps. This provides an
isotopy $s_k'{}^t$ between $s_k^0$ and $s_k^1$. The foliations
$\a_k^t=(\q_k'{}^t)^*\a$, $\q_k'{}^t=\PP(s'_k{}^t)$, $t\in [0,1]$,
provide an isotopy between the initial ones, as required.
\end{proof}

\begin{remark}
 The perturbation of the foliation carried out in section
 \ref{sec:ah-foli} can be done in a one-parameter family $\a_k^t$,
 as long as we start with a one-parameter family of asymptotically
 holomorphic functions $f_k^t$ and $h_k^t$ to start with.
 Therefore for a family of $(\g,\e)$-regular asymptotically
 holomorphic foliations $\a_k^t$, there exists a family of
 symplectic foliations $\hat{\a}_k^t$ interpolating between the
 perturbations $\hat{\a}_k^0$, $\hat{\a}^1_k$ of $\a_k^0$,
 $\a_k^1$ carried out in proposition \ref{prop:holo_sim}. So the
 construction of the symplectic foliations is unique up to
 symplectic isotopy, for $k$ large enough.
\end{remark}

\begin{remark} \label{rem:non-integer}
 Suppose that $(M,\o)$ is a symplectic manifold with $[\o]/2\pi$
 not an integer cohomology class in $H^2(M;\R)$. Choose a
 compatible almost complex structure $J$. We may take a small
 perturbation $\o'$ of $\o$ which is still symplectic and
 compatible with $J$, such that $[\o']/2\pi$ is a rational
 cohomology class. Therefore there is a positive integer $M$ such
 that $M[\o']/2\pi$ is an integer cohomology class.

 Applying the theorem above for $M\o'$ we get asymptotically
 holomorphic foliations, and therefore symplectic foliations for
 $(M,\o)$, with hermitian normal bundles $L^{\ox Nk}$ where
 $c_1(L)=M [\o']/2\pi$.
\end{remark}

\section{Examples}\label{sec:examples}
 We can apply all the precedent constructions to any fixed
 foliation in the projective space. There is a large number of
 examples \cite{Ca94,CL96}. We are going to compute explicitly
 some classical cases.

\subsection{Application to Lefschetz pencils}
 First, we can recover Donaldson's result \cite{Do99} on the
 existence of Lefschetz pencils. We need the following definition.

\begin{definition}
A branched $(p,q)$ Lefschetz pencil, with $p,q>0$ relatively
prime, over an oriented closed manifold $M$ consists of the
following set of data:
\begin{enumerate}
\item A codimension four smooth submanifold $B$.
\item A map $f:M-B \to \CP^1$ which is a submersion outside a finite
 set of points $\Delta$.
\end{enumerate}
 Also the data fit in the following models
\begin{itemize}
\item Given a point $x\in B$, there exists a neighborhood of $x$ with
 oriented coordinates $(z_1, \ldots, z_n)$ of $M$ where the map $f$
 can be written as $f(z_1, \ldots, z_n)={z_2^q}/{z_1^p}$.
\item Given a point $x\in \Delta$, there exists a neighborhood of $x$
 with oriented coordinates $(z_1, \ldots, z_n)$ of $M$ where the
 map $f$ can be written as $f(z_1, \ldots, z_n)=z_1^2 + \cdots +
 z_n^2$.
\end{itemize}
\end{definition}

A branched $(1,1)$ Lefschetz pencil is called a simple Lefschetz
pencil (or a Lefschetz pencil, for brevity). The main result in
\cite{Do99} is

\begin{theorem}
 Let $(M,\o)$ be a symplectic manifold of integer class. There
 exists an integer $k_0>0$ such that for any $k>k_0$,
 $M$ admits a Lefschetz pencil structure $(f_k, B_k, \Delta_k)$
 where all the fibers of the map $f_k$ are symplectic and Poincare
 dual to $k[\o]/2\pi$ and $B_k$ is also symplectic.
\end{theorem}

Donaldson constructs two asymptotically holomorphic sections
$s_k^1, s_k^2 \in L^{\otimes k}$ satisfying certain transversality
properties. The map $f_k$ is defined by $f_k={s_k^2}/{s_k^1}$, and
for $k$ large, it satisfies the required properties. Moreover, the
form $\alpha_k= s_k^1 \, d s_k^2 - s_k^2 \, ds_k^1$ is an
asymptotically holomorphic symplectic foliation. This follows from
the fact that $\alpha_k$ is just a rescaling of the differential
of $f_k$ which is obviously defining a symplectic foliation
whenever it is well defined. The rescaling is performed in order
to $\a_k$ be defined all over the manifold. Moreover the foliation
is as well symplectic in $B_k$, and so it is symplectic all over
$M$. We also can prove

\begin{theorem}
 Let $(M,\o)$ be a symplectic manifold of integer class. There
 exists an integer $k_0>0$ such that for any $k>k_0$, $M$ admits a
 branched $(p,q)$ Lefschetz pencil structure $(f_k, B_k, \Delta_k)$
 where all the fibers of the map $f_k$ are symplectic and Poincare
 dual to $k(p+q)[\o]/2\pi$, and $B_k$ is also symplectic.
\end{theorem}

\begin{proof}
 Fix sections $s_1 \in H^0(T^*\CP^d \ox H^{\ox q})$ and $s_2 \in
 H^0(T^*\CP^d\ox H^{\ox p})$, where $H$ is the hyperplane line
 bundle over $\CP^{d}$. Moreover, suppose that these sections are
 transverse to zero and $Z(s_1) \cap Z(s_2)$ is a transverse
 intersection. Therefore
 $$
  f= \frac{s_{2}^{\ox q}} {s_{1}^{\ox p}}
 $$
 defines a branched $(p,q)$ Lefschetz pencil over $\CP^{d}$.
 Moreover
 $$
 \a=qs_1\,ds_2-ps_2\,ds_1 \in H^0(T^*\CP^d\ox H^{\ox (p+q)})
 $$
is a holomorphic foliation on $\CP^{d}$ satisfying the hypothesis
required in section \ref{sec:existence}. Then by theorem
\ref{thm:main2}, there exists an embedding $\p_k$ of $M$ in
$\CP^{d}$ such that $\a_k=\p_k^* \a$ is an asymptotically
holomorphic and $(\g,\f)$-regular foliation, for uniform
$\g,\f>0$. The perturbation performed in proposition
\ref{prop:holo_sim} takes place well away from the singular locus
and changes $f_k=\p_k^* f$ into an integrating function with a
suitable form around the critical points. Therefore this
perturbation may be done by perturbing either $s_k^1=\p_k^* s_1$
or $s_k^2=\p_k^* s_2$ (since one of them is non-zero). We obtain a
symplectic foliation and the map $f_k=(s_k^2)^{\ox q}/(s_k^1)^{\ox
p}$ defines a branched $(p,q)$ Lefschetz pencil.
\end{proof}

\subsection{Deformations of Lefschetz pencils with non-trivial
holonomy} We may deform the Lefschetz pencils in the presence of
fundamental group as in the algebraic case. We say that a
symplectic foliation $\a$ on a symplectic manifold $(M,\o)$ is a
deformed Lefschetz pencil if there is a connected smooth
codimension four symplectic submanifold $B\subset M$ such that
\begin{enumerate}
 \item Given a point $x\in B$, there are adapted coordinates
 $(z_1, \ldots, z_n)$ around $x$ where the leaves of the foliation
 are of the form $z_2/z_1=\l$.
 \item There is a finite set of critical points $x_j\in M-B$ such that
 at any $x_j$ there are adapted coordinates
 $(z_1, \ldots, z_n)$ where the leaves of the foliation are of the
 form $z_1^2 + \cdots + z_n^2=\l$.
\end{enumerate}

Suppose $\a$ is a deformed Lefschetz pencil with base locus $B$.
The holonomy $H:\pi_1(B)\to \operatorname{PU}(2)$ is defined as
follows. Fix $p_0\in B$ and consider a small transversal
$2$-dimensional disk $\Delta$ to $B$. Identify $\PP(\D)=\CP^1$.
For any loop $\varsigma$ and any $\l\in \CP^1$, lift the path
$\varsigma$ to a path in a tubular neighborhood of $B$ inside the
leaf of $\a$ corresponding to the value $\l$. The endpoint is
defined to be $H(\varsigma)(\l) \in \CP^1$.

\begin{theorem}
Let $(M,\o)$ be a symplectic manifold of integer class such that
$\dim M=2n\geq 6$ and $\pi_1(M)\neq 1$. Let $L\to M$ be the
complex line bundle with $c_1(L)=[\o]/2\pi$. Then for $k$ large
enough there are deformed Lefschetz pencils $\a_k\in
\SC^{\infty}(T^*M\ox L^{\ox k})$ with non-trivial holonomy.
\end{theorem}

\begin{proof}
Let $\rho:\pi_1(M)\to SU(2)$ be a representation and $E_\rho \to
M$ the corresponding flat $\C^2$-bundle. Suppose that $\rho$ is a
small deformation of the trivial representation, so that $E_\rho$
is a topologically trivial bundle. We may understand $E_\rho =
M\x\C^2$ with a flat connection $\nabla_\rho=\nabla +\varpi$,
where $\varpi \in \O^1(\mathfrak{su}(2))$. Now
$|\nabla_\rho-\nabla|_{g_k}= |\varpi|_{g_k}= k^{-1/2}|\varpi|_g$,
i.e., $\nabla_\rho$ is ``asymptotically trivial'' connection.

Endow $M$ with a compatible almost complex structure and let $L\to
M$ be the hermitian line bundle with connection with curvature
$-i\o$. As in the previous section, there are asymptotically
holomorphic sections $s_k=(s_k^1,s_k^2)$ of $\C^2\ox L^{\ox k}$
such that $f_k=s_k^2/s_k^1$ is a symplectic Lefschetz pencil for
$k$ large enough. Let $\a_k$ be the associated foliation with base
locus $B_k=Z(s_k)$. Consider the morphism
 $$
   M-B_k \stackrel{(1,f_k)}{\longrightarrow} M \x \CP^1 \cong \PP
   (E_\rho).
 $$
We pull back the flat distribution of $\PP(E_\rho)$ under this map
to get a foliation $\a_k'$ in $M$.

Let $K_{\a_k}=B_k$ be the Kupka set of $\a_k$ and let $B_\g^k$ be
the neighborhood of $g_k$-radius $\g$ of $K_{\a_k}$. Then in
$M-B_\g^k$ the $1$-form $(\a_k)_{1,0}$ is $\e$-transverse to zero.
On the other hand, the horizontal distribution of $M\x \CP^1$ is
given by $\a=d\l$, where $\l$ is the coordinate in the
$\CP^1$-direction, and the horizontal distribution for
$\PP(E_\rho)$ is given by $\a'=d\l + O(k^{-1/2})$. Since
$\a_k=(1,f_k)^* \a$ and $\a'_k=(1,f_k)^*\a'$ (up to a factor
$(s_k^1)^{\ox 2}$, which is uniformly bounded) we have that
$|\a_k-\a'_k|_{C^1,g_k}=O(k^{-1/2})$. So $\a'_k$ is an
asymptotically holomorphic foliation such that $(\a'_k)_{1,0}$ is
$\e/2$-transverse to zero in $M-B_\g^k$. After the perturbation in
proposition \ref{prop:holo_sim}, it defines a symplectic
foliation.

Now we look at $B_\g^k$, where $f_k$ has transversal type
$z_1dz_2-z_2dz_1$. Let us see that this is stable for small
perturbations. We construct $E_\rho$ in a different way: take the
universal covering space $\pi: \tilde{M} \to M$ and consider
$\tilde{M}\x\C^2$ with the trivial connection. Identify $(x,f)$
with $(h(x), \rho(h) f)$ for any deck transformation $h$, to get
$E_\rho \to M$. Consider the identification of $E_\rho$ with the
trivial bundle as a (not connection preserving) isomorphism $\psi:
\underline{\C}^2 \to E_\rho$. This lifts to $\tilde{\psi}
:\tilde{M}\x\C^2 \to \tilde{M}\x\C^2$, where $\tilde{\psi}^*
\nabla=\nabla +\tilde{\varpi}$, $\tilde{\varpi}=\pi^*\varpi$.
Therefore $\tilde{\varpi}=\tilde{\psi}^{-1} d\tilde{\psi}$ and we
may consider $\tilde{\psi}$ as a map $\tilde{M}\to
\operatorname{SU}(2)$ satisfying $\tilde{\psi}(h(x))= \rho(h)
\tilde{\psi}(x)$. Look at the map
  $$
   \tilde{M}-\tilde{B}_k
   \stackrel{(1,\tilde{f}_k)}{\longrightarrow} \tilde{M} \x \CP^1
   \stackrel{\tilde{\psi}}{\isom} \tilde{M} \x \CP^1,
  $$
where $\tilde{f}_k=f_k\circ \pi$ and $\tilde{B}_k=\pi^{-1}(B_k)$.
Fix a point $p\in B_k$ and $(z_1,z_2,\ldots, z_n)$ coordinates
around $p$ such that $K_{\a_k}=\{z_1=z_2=0\}$ and $f_k=z_2/z_1$.
Looking at any point in $\tilde{M}$ over $p$ we see that the
leaves of $\a'_k$ are the level sets of $\tilde{\psi}\circ
\tilde{f}_k$. Denoting
  $$
  {w_1 \choose w_2}=\tilde{\psi}(z_1,\ldots, z_n)
  {z_1 \choose z_2},
  $$
we have new adapted coordinates $(w_1,w_2,z_3,\ldots, z_n)$ such
that the leaves of $\a'_k$ are of the form $w_2/w_1=\l$. Note that
the change of coordinates is asymptotically holomorphic, since
$|d\tilde{\psi}|_{C^1}=O(k^{-1/2})$.

We have a symplectic foliation $\a'_k$ with $K_{\a_k}=B_k$ a
symplectic smooth submanifold of codimension $4$. Let us see that
this foliation has non-trivial holonomy (and therefore it is not a
Lefschetz pencil). By the Lefschetz theorem in \cite{Au97}, $B_k$
is connected and $\pi_1(B_k)\surj \pi_1(M)$ for $k$ large enough.
Therefore $\rho$ defines a non-trivial representation of
$\pi_1(B_k)$ which determines $\PP(E_\rho)|_{B_k}$. Fix $p\in B_k$
and $\D$ a small transversal $2$-dimensional disk to $B_k$ at $p$.
Identify $\PP(\D)=\CP^1$. Let $\varsigma \in \pi_1(B_k)$ be a loop
and $\l \in \CP^1$. We move $\varsigma$ into a regular leaf
starting at the leaf determined by $\l$ at $p$. Then looking at
the picture in $\tilde{M}$ we see that the end-point is
$\rho(\varsigma) (\l)$. Note that if $\rho(\pi_1(M))\l\subset
\CP^1$ is infinite then the leaf corresponding to $\l$ is not
compact. Moreover if $\rho$ is not the identity in
$\operatorname{PU}(2)$, then (the closures of) the leaves are not
smooth at $B_k$ and therefore $\a'_k$ does not define a Lefschetz
pencil.
\end{proof}

\subsection{Application to asymptotically holomorphic logarithmic foliations}
Another useful example of foliations in the projective space are
the logarithmic foliations as discussed in subsection
\ref{subsec:log}.

\begin{definition}
Let $M$ be a closed manifold, $L_1,\ldots, L_p$ a family of
complex line bundles over $M$ and $\lambda_1, \ldots, \lambda_p\in
\C$ complex numbers such that $\sum \l_i c_1(L_i) =0 \in
H^2(M;\C)$. Choose sections $f_1, \ldots, f_p$ of the bundles
$L_1, \ldots, L_p$. Then a logarithmic foliation with normal
bundle $L=L_1 \otimes \cdots \otimes L_p$ is given by the twisted
$1$-form
 $$
  \a=f_1\cdots f_p \sum_{i=1}^p \l_i \frac{df_i}{f_i}
  \in \SC^\infty(M,T^*M\ox L).
 $$
\end{definition}

The condition above ensures that in the open subset $\{f_1\neq 0,
\ldots, f_p\neq 0\}$ we have $\a=f_1\cdots f_p \, d(\log F)$,
where $F=f_1^{\l_1}\cdots f_p^{\l_p}$ is a function. It is easy to
see that this gives a well defined $L$-valued $1$-form all over
$M$. Now consider the manifold and the sections to be holomorphic.
We say that the foliation is generic if:
\begin{enumerate}
 \item $p\geq 3$ and for any $i=1, \ldots, p$, the line bundle
  $L_i$ is positive.
 \item For any $i=1, \ldots,p$, the hypersurface defined
  by the equation $\{f_i=0\}$ is irreducible and
  $\{f_1 \cdots f_p=0 \}$ is a divisor with normal crossings.
\end{enumerate}

In the particular case of the projective space we have that $f_i
\in H^0(\CP^d, \SO(n_i))$ is a homogeneous polynomial of degree
$n_i$. The condition $\sum \l_i c_1(L_i) =0$ translates into
 \begin{equation}
  \sum_{i=1}^p n_i\l_i=0. \label{cuadrop}
 \end{equation}
With this kind of foliations at hand we can prove

\begin{theorem}
Let $(M,\o)$ be a symplectic manifold of integer class and let
$L\to M$ be a complex line bundle with $c_1(L)=[\o]/2\pi$. Given
$(n_1, \ldots, n_p)$ positive integers and $(\l_1, \ldots, \l_p)$
satisfying the condition \eqref{cuadrop}, then for $k$ large
enough there exists an asymptotically holomorphic sequence of
sections $(f_k^1, \ldots, f_k^p)$ of $L^{\ox kn_1}\oplus \cdots
\oplus L^{\ox kn_p}$ such that the associated logarithmic
foliation
 $$
 \a_k=f_k^1\cdots f_k^p\sum_{i=1}^p \l_i\frac{df_k^i}{f_k^i}
 \in \SC^\infty(M, T^*M \ox L^{\ox kN}),
 $$
where $N= n_1+ \cdots +n_p$, is a symplectic foliation.
\end{theorem}

\begin{proof}
Choose a generic logarithmic foliation $\a$ in the projective
space with polynomials $(f_1,\ldots, f_p)$ of degrees $(n_1,
\ldots, n_p)$ and complex numbers $(\lambda_1, \ldots, \lambda_p)$
as in the statement. Then the foliation satisfies the conditions
of theorem \ref{thm:main2}, so there is a family of asymptotically
holomorphic embeddings in the projective space $\phi_k:M \to
\CP^d$ such that the pull-back foliations $\phi_k^* \alpha$ give
an asymptotically holomorphic sequence of foliations which are
$(\g,\f)$-regular. The perturbation of proposition
\ref{prop:holo_sim} moves the integrating factor $\log
\left((f_k^1)^{\l_1}\cdots (f_k^p)^{\l_p}\right)$, where
$f_k^i=\phi_k^* f_i$. When the critical point is well away from
every $D_k^i=Z(f_k^i)$ this perturbation can be absorbed into a
perturbation of some $f_k^i$. To avoid that the critical points
get close to $D_k^i$, just take the embeddings $\p_k$ to be
transverse to every $D_i=\{ f_i=0\}\subset \CP^d$ by using theorem
\ref{good_inter1}. This produces logarithmic symplectic foliations
for $k$ large enough.
\end{proof}

Remark that any element of this family of foliations is not
equivalent to a Lefschetz pencil. In particular, they are not of
Kupka type.

\end{document}